\documentclass[12pt]{article}

\usepackage{a4wide}
\usepackage{amssymb,amsmath}

\allowdisplaybreaks[1]

\newcommand{\titre}{Quantum unique factorisation domains}

\newenvironment{proof}{\begin{trivlist}\item[]{\it
Proof.}}{\hfill$\square$\end{trivlist}}

\newtheorem{theorem}{Theorem}[section]
\newtheorem{corollary}[theorem]{Corollary}
\newtheorem{definition}[theorem]{Definition}
\newtheorem{hypothesis}[theorem]{Hypothesis}
\newtheorem{lemma}[theorem]{Lemma}
\newtheorem{proposition}[theorem]{Proposition}
\newtheorem{remark}[theorem]{Remark}

\def\sse{\subseteq}
\def\goesto{\longrightarrow}

\def\mc{{\mathbb{C}}}

\def\mz{{\mathbb{Z}}}
\def\k{k}
\def\w{w}

\def\wr{\widehat{R}}

\def\detq{{\rm det}_q}

\def\ch{{\mathcal H}}
\def\co{{\mathcal O}}

\def\gq{{\cal G}_q}
\def\gqmn{\gq(m,n)}

\def\oq{{\cal O}_q}
\def\oqmn{\co_q(M_n)}
\def\oqmmn{\co_q(M_{m,n})}
\def\oqmmn{\co_q(M_{m,n})}
\def\oqmnm{\co_q(M_{n,m})}

\begin{document}

\title{\titre}
\author{S Launois, T H Lenagan and L
Rigal
\thanks{This work 
was supported by Leverhulme Research Interchange
Grant F/00158/X }\;
}
\date{}

\maketitle


\begin{abstract}
We prove a general theorem showing that iterated skew polynomial extensions of the type
which fit the conditions needed by Cauchon's deleting derivations theory and
by the Goodearl-Letzter stratification theory are unique factorisation rings
in the sense of Chatters and Jordan. This general result applies to many
quantum algebras; in particular, generic quantum matrices and quantized
enveloping algebras of the nilpotent part of a semisimple Lie algebra are
unique factorisation domains in the sense of Chatters. 
By using noncommutative
dehomogenisation, the result also extends to generic quantum grassmannians.

\end{abstract}



\section*{Introduction}\label{sec:intro}

In \cite{ch}, Chatters introduced the notion of a noncommutative unique
factorisation domain in the following way. An element $p$ of a noetherian
domain $R$ is said to be {\em prime} if (i) $pR = Rp$, (ii) $pR$ is a height
one prime ideal of $R$, and (iii) $R/pR$ is an integral domain. A noetherian
domain $R$ is then said to be a {\em unique factorisation domain}, noetherian UFD for
short, if $R$ has at least one height one prime ideal, and every height one
prime ideal is generated by a prime element. As well as the usual commutative
noetherian UFDs, examples include universal enveloping algebras of finite dimensional
solvable Lie algebras over $\mc$. However, one of the deficiences of this
definition is that the class of noetherian UFDs is not closed under polynomial
extensions, as \cite[Example 2.11]{ch} shows. The problem is that the
condition of height one prime factors being domains does not pass up to
polynomial extensions.

In order to remedy this deficiency, in a later paper, the notion of a {\em
noetherian unique factorisation ring}, noetherian UFR for short, was introduced by
Chatters and Jordan, \cite{cj}. For a large class of rings (namely, the
noetherian prime rings satisfying the descending chain condition on prime
ideals), being a noetherian unique factorisation ring amounts having height
one primes principal (that is generated by a normal element). This condition
is closed under polynomial extensions, and, indeed, they then are able to
prove theorems about skew polynomial extensions of the type $R[x;\sigma]$ and
$R[x;\delta]$. However, they do not prove any results about general skew
polynomial extensions of type $R[x;\sigma, \delta]$.

In many quantum algebras, in the generic case where the deformation parameter
$q$ is not a root of unity, it is known that all prime ideals are completely
prime, and then the distinction between a noetherian domain being a noetherian UFD and a
noetherian UFR disappears and so the results of \cite{cj} on noetherian UFRs also apply to noetherian UFDs in
this setting.

The purpose of this paper is to obtain a theorem on unique factorisation for
certain extensions of the type $R[x;\sigma, \delta]$ that arise naturally in
the study of quantum algebras. Once this theorem is proved, an iterated
version is obtained which is sufficient to show that many quantum algebras are
noetherian UFDs. In particular, we show that the algebra of generic quantum matrices,
$\oqmmn$ is a noetherian UFD, as is the quantized enveloping algebra $U_q^+(g)$.

Roughly speaking, an iterated skew polynomial extension will be a noetherian UFD provided 
that the Cauchon theory of deleting derivations can be applied, 
and that there is a torus action $\ch$ for which the Goodearl-Letzer 
stratification theory applies. Exact requirements will be given as 
they become necessary.

In the case of quantum matrices, we can go further, since we can identify 
the height one prime ideals that are $\ch$-primes for the natural torus 
that acts. 

In the final section, we use the idea of noncommutative dehomogenisation, 
developed in \cite{klr} to deduce that generic quantum grassmannians are 
noetherian UFDs.   

For general results concerning noetherian rings and localisation, we refer the
reader to \cite{gw} or \cite{mr}.

Throughout the paper, $k$ denotes a field. 

\section{Non commutative unique factorisation rings}\label{section:ufr}

This section investigates the behaviour
of the notion of a {\it noetherian unique factorisation ring}, as defined
in \cite{cj} by Chatters and Jordan, under localisation by normal elements.\\

To start with, we recall the definition of noetherian unique factorisation ring; 
further details concerning this notion can be found in \cite{cj}. \\

An ideal $I$ in a ring 
$A$ is called {\it principal} if there exists a normal element $x$ in 
$A$ such that $I=\langle x \rangle \; (= xA = Ax)$.

\begin{definition} \label{def-UFR}
A ring $A$ is called a noetherian unique factorisation ring (noetherian UFR for short) 
if: 
(i) $A$ is a prime noetherian ring, and \\
(ii) any nonzero prime ideal in $A$ contains a nonzero principal prime 
ideal.
\end{definition}

\begin{definition} \label{def-UFD}
A noetherian UFR $A$ is said to be a unique factorisation domain (noetherian UFD for short) if $A$
is a domain and each height one prime ideal $P$ of $A$ is completely prime;
that is, $A/P$ is a domain for each  height one prime ideal $P$ of $A$. 
\end{definition}

\begin{remark} \rm 
If $A$ is a prime noetherian ring that satisfies the descending chain
condition for prime ideals, then $A$ is a noetherian UFR if and only if height
one primes are principal (see \cite{cj}). Hence, the notions of noetherian UFR
and noetherian UFD are good generalisations of the usual notion of unique factorisation
domain for commutative rings (see in particular Corollaries 10.3 and 10.6 in
\cite{e}).

Note that the algebras  
we are dealing with are all noetherian and have 
finite Gelfand-Kirillov dimension; so, they satisfy the descending chain
condition for prime ideals, see for example, \cite[Corollary 3.16]{kl}.
\end{remark}

We start by proving a noncommutative analogue of Nagata's Lemma (in the
commutative case, see \cite{e} 19.20 p.\! 487). The following result is taken
from \cite{dr}, where it appears without proof. We include a proof here, for
the convenience of the reader, since it is crucial to a part of our argument.

If $A$ is a prime noetherian ring and $x$ a nonzero normal element 
of $A$, we denote by $A_x$ the right localisation of $A$ with respect to the 
powers of $x$. 

\begin{lemma} \label{nagata}
Let $A$ be a prime noetherian ring and $x$ a nonzero, nonunit, normal element 
of $A$ such that $\langle x \rangle$ is a completely prime ideal of $A$.\\
(i) If $P$ is a prime ideal of $A$ not containing $x$ and such that the prime 
ideal $PA_{x}$ of $A_{x}$ is principal, then $P$ is principal.\\
(ii) If $A_{x}$ is a noetherian unique factorisation ring, then so is $A$. \\
(iii) If  $A_{x}$ is a noetherian unique factorisation domain, then so is $A$.
\\
\end{lemma}

\begin{proof}~~
(i) The result is trivial if $P=0$, so we assume that $P\neq 0$. Since $x$ is
a nonzero normal element of the prime ring $A$ one may localise $A$ with
respect to the multiplicative set of powers of $x$ and there is canonical
embedding $A \hookrightarrow A_{x}$. Moreover, $Q := PA_{x}$ is a prime ideal
of $A_{x}$ whose contraction to $A$ is $P$, since $P$ is a prime ideal of $A$
not containing $x$ . Let us suppose that $Q$ is a principal ideal. Then,
clearly, there exists $q \in A$, normal in $A_{x}$, such that $Q=qA_{x}$.
Moreover, one may assume the right ideal $qA$ maximal for this property, since
$A$ is right noetherian. Suppose that $q\in Ax$. Then there exists $p$ in $A$
such that $q=px$ (in particular $qA \subseteq pA$). But then, $Q=pA_{x}$ and
$p$ is normal in $A_{x}$. The maximality of $qA$ leads to $qA=pA$ from which
follows the existence of $r \in A$ such that $p=qr$ and hence $q=qrx$. Since
$q$ is a non-zero normal element in the prime ring $A_{x}$, the above equality
gives $1=rx$ (with $r \in A$), a contradiction, since $x$ is not a unit. Thus,
$q \notin Ax$. Now, let $p \in P \subseteq Q$; so that there exist $r \in A$
and $t \in {\mathbb N}$ with $p=qrx^{-t}$, and we may choose $t$ minimal for
this property. If $t>0$ then $r \notin Ax$, by the minimality of $t$. The
above equation then leads to $px^{t}=qr$; and so either $q$ or $r$ must be in
$Ax$ which is a contradiction. Thus, $t=0$ and so $p \in qA$. Hence, $P
\subseteq qA$. Also, $qA \subseteq qA_{x} \cap A =Q \cap A =P$; so that $P =
qA$. A similar argument gives $P=Aq$. Hence $P = Aq = qA$ which proves the
first claim.\\ 
(ii) Let us now assume that $A_{x}$ is a noetherian UFR. If
$Q_{0}$ is a non-zero prime ideal of $A$ not containing $x$, then $Q_{0}A_{x}$
is a non-zero prime ideal of $A_{x}$. Since $A_{x}$ is a noetherian UFR,
$Q_{0}A_{x}$ contains a nonzero principal prime ideal $P$ which is the
extension to $A_{x}$ of its contraction $P_{0}$ in $A$. By part (i), the ideal
$P_0$ is principal, since $P$ is principal. Thus, $P_0$ is a nonzero principal
prime ideal contained in $Q_0$. Moreover, if $Q_{0}$ is a prime ideal of $A$
containing $x$, then it contains the nonzero principal prime ideal $\langle x
\rangle$. We have proved that each nonzero prime ideal of $A$ contains a
nonzero principal prime ideal, which means that $A$ is a noetherian UFR. 
\\
(iii)
Suppose that $A_x$ is a noetherian UFD. Then part (ii) shows that $A$ is a noetherian UFR. Let $P$
be a prime ideal of height one in $A$. If $x\in P$ then $P = \langle
x\rangle$ and so $P$ is completely prime, by assumption. 
Otherwise, standard localisation theory shows that $PA_x$ is a
prime ideal of height one in $A_x$ and that $P = PA_x \cap A$. Thus, $A/P$
embeds in $A_x/PA_x$, which is a domain; and so $A/P$ is a domain, as
required.
\end{proof}

Proposition \ref{pull-back-UFR} below will be of central use later. It gives a
way to pull back the unique factorisation property from a certain type of
localisation to the initial ring. The following lemma is needed in the proof
of the proposition.

\begin{lemma} \label{pre-pull-back-UFR}
Let $R$ be a prime noetherian ring and suppose that $d,s$ are normal
elements of $R$ such that $dR$ is prime and $s\notin dR$. Then, there exist
units $u,v \in R$ such that $ds=sdu$ and $sd=vds$.
\end{lemma}

\begin{proof}~~ 
If either $d$ or $s$ is zero, then the result is trivial; so we assume that 
$d,s
\neq 0$. Since $s$ is normal in a prime ring, $s$ is regular and we can 
associate to it an automorphism $\sigma \, : \, R \longrightarrow R$ 
such that $xs=s\sigma(x)$, for all $x\in R$. 
Set $P:=dR = Rd$. Then 
$s\sigma(P) = Ps \subseteq P$; and so $\sigma(P) \subseteq P$, since $s$ 
is normal and not in $P$. 
Hence, $P \subseteq \sigma^{-1}(P)$, and it follows
that there is an ascending chain $P \subseteq \sigma^{-1}(P) \subseteq 
\sigma^{-2}(P) \subseteq \dots$ of ideals of $R$. 
The noetherian hypothesis then ensures that
there exists $n \in {\mathbb N}$ such that $\sigma^{-n}(P)=
\sigma^{-(n+1)}(P)$, and so $\sigma(P)=P$; that is,  
$\sigma(d)R=dR$. From this it follows that $dsR=sdR$, which
gives the existence of $u,u'\in R$ 
such that $ds=sdu$ and $sd=dsu'$. But then,
$ds= sdu = dsu'u$; and so $u'u =1$ which shows $u$ is a unit in $R$.
We also have $Rds=Rsd$, since $d$ and $s$ are normal, and it follows 
in a similar manner that there exists
a unit $v$ in $R$ such that $sd=vds$.
\end{proof}

\begin{proposition} \label{pull-back-UFR}
Let $R$ be a prime noetherian ring and suppose that $d_1, \dots, d_t$ are
nonzero normal elements of $R$ such that the ideals $d_1R,\dots,d_tR$ are
completely prime and pairwise distinct. Denote by $T$ the right quotient ring
of $R$ with respect to the right denominator set generated by $d_1,\dots,d_t$.
If $T$ is a noetherian UFR then so is $R$. Also, if $T$ is a noetherian UFD then so is
$R$. 
\end{proposition}

\begin{proof}~~ 
We proceed by induction on $t$, the result being true for $t=1$ by
Lemma~\ref{nagata}~(ii).\\
Assume that the result is true up to order $t \in {\mathbb N}^\ast$.
We will work in the right quotient ring of fractions of $R$ in which all the
localisations of $R$ are naturally embedded.
Denote by ${\mathcal S}_{t+1}$ the multiplicative subset of $R$ generated by
$d_1,\dots,d_{t+1}$ and by ${\mathcal S}_t$ the multiplicative subset of $R$
generated by $d_1,\dots,d_t$. Hence $T=R{\mathcal S}_{t+1}^{-1}$. 
We first show, using the above lemma, that $d_{t+1}$ is a nonzero normal 
element of $R{\mathcal S}_t^{-1}$. Let $(a,s)\in R\times{\mathcal S}_t$;
hence $s$ is normal in $R$ and, due to the 
hypothesis that the ideals $d_iR$ are
completely prime and pairwise distinct, $s\notin d_{t+1}R$ (by the principal
ideal theorem). So, by the lemma above, there exist elements $u,v \in R$ such 
that $d_{t+1}s=sd_{t+1}u$ and $sd_{t+1}=vd_{t+1}s$. In addition, since 
$d_{t+1}$ 
is normal in $R$, there exist $b,c\in R$ such that $ad_{t+1}=d_{t+1}b$ and 
$d_{t+1}a=cd_{t+1}$. Hence, we have
$as^{-1}d_{t+1}=ad_{t+1}us^{-1}=d_{t+1}bus^{-1}$ and 
$d_{t+1}as^{-1}=cd_{t+1}s^{-1}=cs^{-1}vd_{t+1}$. 
It follows that $d_{t+1}$ is indeed a nonzero normal element of 
$R{\mathcal S}_t^{-1}$.

Let ${\mathcal S}$ be the multiplicative subset of $R{\mathcal S}_t^{-1}$
generated by $d_{t+1}$. Notice that, $R{\mathcal S}_{t+1}^{-1} 
= (R{\mathcal S}_t^{-1}){\mathcal S}^{-1}$, as is easily verified. 
Of course, $R{\mathcal S}_t^{-1}$ is prime noetherian and 
the ideal $d_{t+1}R{\mathcal S}_t^{-1}$ is completely prime since $d_{t+1}R$ is 
completely prime and does not intersect ${\mathcal S}_t$.

Now assume that $T=(R{\mathcal S}_t^{-1}){\mathcal S}^{-1}$ is a noetherian
UFR. By the comments above, Lemma~\ref{nagata}~(ii) may be applied and we get
that $R{\mathcal S}_t^{-1}$ is a noetherian UFR. Now, the induction hypothesis
gives that $R$ is a noetherian UFR, as required.

Finally, suppose that $T$ is a noetherian UFD. Then $T$ is certainly a
noetherian UFR; and so $R$ is a noetherian UFR, by the first part of this
result. That $R$ is a noetherian UFD then follows by standard localisation
theory (cf. the proof of Lemma~\ref{nagata}~(iii)).
\end{proof}

\section{Height one $\ch$-primes in Cauchon extensions} 

Most of the algebras that we are considering in this paper have groups acting
on them in natural ways. The study of the prime spectra of such algebras is
often facilitated by first studying ideals invariant under the natural group
action. We begin this section by recalling some standard terminology
concerning ideals invariant under group actions. A convenient reference is
\cite[II.1.8, II.1.9]{bg}. Let $\ch$ be a group acting by automorphisms on a
ring $R$. An ideal $I$ of $R$ is an {\em $\ch$-ideal} 
provided that $h(I) = I$ for
all $h\in \ch$. A proper $\ch$-ideal is an 
{\em $\ch$-prime ideal} provided that whenever
$IJ\sse P$ for $\ch$-ideals $I,J$ of $R$ then either $I\sse P$ or $J\sse P$.
The set of $\ch$-prime ideals of $R$ is denoted by 
{\em $\ch{\rm -Spec}(R)$}. It is obvious
that a prime ideal $P$ that is an $\ch$-ideal is an $\ch$-prime ideal. The
converse is not true in general; however, it will usually be true for the
algebras that interest us in this paper (see comments after
Definition~\ref{def-cgl}).

\begin{hypothesis}\label{cauchon1} 
Let $A$ be a domain that is a noetherian $k$-algebra and suppose that $\sigma$
is a $k$-automorphism of $A$. Suppose that there is a group $\ch$ acting as
automorphisms on the skew Laurent extension $A[X^{\pm 1}; \sigma]$ in such a
way that $X$ is an $\ch$-eigenvector and $A$ is stable under $\ch$. Further,
suppose that the action of $\sigma$ on $A$ coincides with the action of an
element $h_0 \in \ch$. Finally, suppose that there is a non root of unity
$\lambda_0$ in $k^*$ such that $h_0.X = \lambda_0 X$. 
\end{hypothesis}

Given the conditions of this hypothesis, we are going to show that there is a
bijection between the $\ch$-ideals of $A$ and the $\ch$-ideals of $A[X^{\pm 1};
\sigma]$, and, consequently, there is a bijection between $\ch{\rm -Spec}(A)$
and those $\ch$-primes of $A[X;\sigma]$ that do not contain $X$.

\begin{lemma}
\label{lem21} 
Assume Hypothesis~\ref{cauchon1} 
and let $I$  be an $\ch$-invariant ideal of $A[X^{\pm 1}; \sigma]$.
Suppose that $x = a_1X^{k_1} +\dots
+a_nX^{k_n}\in I$, with $a_i\in A$ and $k_i$ all distinct. Then, each $a_i \in
I \cap A$.  Consequently, 
$I = \oplus_{i\in \mz}\, (I\cap A)X^i$.
\end{lemma}

\begin{proof} 
The proof is by induction on $n$. If $n=1$ the result is trivial, since $X$ is
invertible. Suppose now that $n>1$. Since $I$ is an $\ch$-ideal, the
element $Xx-\lambda_0^{- k_n}h_0(x) X$ belongs to $I$. However,
\begin{eqnarray*}
Xx-\lambda_0^{- k_n}h_0(x) X &=& 
\sum_{i=1}^{n} h_0(a_i) X^{k_i+1}-\sum_{i=1}^{n}
\lambda_0^{k_i-k_n} h_0(a_i) X^{k_i+1}\\
&=& 
\sum_{i=1}^{n-1}
(1-\lambda_0^{(k_i-k_n)})h_0(a_i) X^{k_i+1};
\end{eqnarray*}
so that $\sum_{i=1}^{n-1}
(1-\lambda_0^{(k_i-k_n)})h_0(a_i) X^{k_i+1} \in I$.
By the induction hypothesis, we see that $(1-\lambda_0^{(k_i-k_n)})h_0(a_i) \in I$ for each $1\leq i\leq
n-1$. The elements $(1-\lambda_0^{(k_i-k_n)})$ are nonzero, since $\lambda_0$
is not a root of unity and the $k_i$ are distinct. Thus, each $h_0(a_i)$ is in
the $\ch$-ideal $I\cap A$, and so each $a_i \in I \cap A$ for $1\leq
i\leq n-1$.
Finally,
$a_n X^{k_n}=x -a_1 X^{k_1}-\dots-a_{n-1} X^{k_{n-1}} \in I$; and so $a_n \in
I\cap A$ also.
\end{proof}

The next result follows easily from this lemma.

\begin{theorem} 
Assume  Hypothesis~\ref{cauchon1}. 
Then there is an inclusion preserving bijection from the set of
$\ch$-ideals of $A$ 
to the set of $\ch$-ideals 
of $A[X^{\pm 1};\sigma]$ given by 
$I \mapsto \oplus_{i\in\mz}\, IX^i$; 
its inverse is defined by $J \mapsto J\cap A$.
Furthermore, these bijections induce order preserving bijections between
$\ch{\rm -Spec}(A)$ and $\ch{\rm -Spec}A[X^{\pm 1};
\sigma]$. 
\end{theorem}

Let $\ch$ be a group acting by automorphisms on a noetherian 
ring $R$ and suppose that
$X$ is a normal $\ch$-eigenvector. Then there is a bijective correspondence
between the $\ch$-prime ideals of $R$ that do not contain $X$ and the
$\ch$-prime ideals of $R[X^{-1}]$, cf \cite[Exercise II.1.J]{bg}. Using this
fact, the next corollary follows easily.

\begin{corollary} \label{cont-ext-cor}
Assume  Hypothesis~\ref{cauchon1}. Then 
contraction $P \mapsto P\cap A$ and extension $P\mapsto \oplus_{i\geq
0}\, PX^i$ provide inverse order preserving bijections between the
$\ch$-prime ideals of $A[X;\sigma]$ that do not contain $X$ and $\ch{\rm
-Spec}(A)$.  
\end{corollary}

\begin{definition}
{\rm 
Let $A$ be a domain that is a noetherian $k$-algebra and let
$R=A[X;\sigma,\delta]$ be a skew polynomial extension of $A$. We say that $R
=A[X;\sigma,\delta]$ is a {\em Cauchon Extension} provided that

\begin{itemize}

\item $\sigma$ is a $k$-algebra automorphism of $A$ and
$\delta$ is a $k$-linear locally nilpotent 
$\sigma$-derivation of $A$. Moreover we assume that there exists 
$q \in k^*$ which is not a root of
unity such that $\sigma \circ \delta = q \delta \circ \sigma$.

\item There exists an abelian group $\ch$ which acts on $R$ by
$k$-algebra 
automorphisms such that $X$ is an $\ch$-eigenvector and $A$ is
$\ch$-stable.

\item $\sigma$ coincides with the action on $A$ of an element $h_0 \in \ch$.

\item Since $X$ is an $\ch$-eigenvector and since $h_0 \in \ch$, there
 exists $\lambda_0 \in k^*$ such that $h_0.X=\lambda_0
 X$. We assume that $\lambda_0$ is not a root of unity.

\item Every $\ch$-prime ideal of $A$ is completely prime. 

\end{itemize}
}
\end{definition} 

Note that the conditions of \cite[II.5.3]{bg} are satisfied by any Cauchon
extension; and so, for example, every $\ch$-prime of $R$ is also completely
prime, by \cite[Proposition II.5.11]{bg}.

In a Cauchon extension $R=A[X;\sigma,\delta]$ 
the set $S=\{X^n \mid n \in \mathbb{N} \}$  is a 
right and left Ore set 
in $R$, \cite[Lemme 2.1]{c1}; and so we can form the Ore localization
$\widehat{R}:=RS^{-1}=S^{-1}R$.

For each $a\in A$, set 

$$\theta (a) = \sum_{n=0}^{+ \infty} \frac{(1-q)^{-n}}{[n]!_q}
\delta^n \circ \sigma^{-n} (a) X^{-n} \in\widehat{R}  $$

Note that $\theta(a)$ is a well-defined element of $\wr$, since
$\delta$ is locally nilpotent, $q$ is not a root of unity, 
and $0\neq 1-q\in k$. 

The following facts are established in \cite[Section 2]{c1}. The map
$\theta:A \goesto \wr$ is a $k$-algebra monomorphism. Let
$A[Y;\sigma]$ be a skew polynomial extension. Then $\theta$ extends to a
monomorphism $\theta: A[Y;\sigma] \goesto \wr$ with $\theta(Y) = X$.
Set $B = \theta(A)$ and $T = \theta(A[Y;\sigma]) \sse \wr$. Then $T =
B[X;\alpha]$, where $\alpha$ is the automorphism of $B$ defined by
$\alpha(\theta(a)) = \theta(\sigma(a))$.

The element $X$ is a normal element in $T$, and so the set $S$ is an
Ore set in $T$ and Cauchon shows that  $TS^{-1}= S^{-1}T = \wr$.

Since $X$ is an $\ch$-eigenvector, it follows from \cite[Exercise
II.1.J]{bg} that $\ch$ also  acts by automorphisms on
$\widehat{R}$. Moreover, the following result shows that 
the group $\ch$ also acts by automorphisms on $T$ and $B$ by
restriction. 

Note, for later use, that, since each element of $B =\theta(A)$ is of the form
$\theta(a) = \sum_{i=0}^n\, a_iX^{-i}$ for some $a_i \in A$, and each element
of $R$ is of the form $\sum_{i=0}^n \, c_iX^i$ for some $c_i \in A$, it
follows that $B\cap R\sse A$.

The next result shows that the action of $\ch$ can be transferred to $B$ via
$\theta$. This result is essentially a generalisation of \cite[Proposition
2.1]{c1}. 

\begin{lemma} Let $R=A[X;\sigma,\delta]$ be a Cauchon extension and let 
$h\in \ch$. Then 
$h.\theta(a)=\theta(h.a)$ for each $a\in A$.

\end{lemma} 

\begin{proof}
 We start by showing inductively that
$h.\delta^n(a)=\lambda_h^n \delta^n(h.a)$ for all $n\in \mathbb{N}$,
$a\in A$ and $h\in \ch$, where $\lambda_h$ denotes the $\ch$-eigenvalue
associated to the $\ch$-eigenvector $X$.

If $n=0$, there is nothing to prove. Now we assume that $n\geq 1$.
Then, since $\delta^n(a)=X \delta^{n-1}(a) - \sigma \circ
\delta^{n-1}(a) X = X \delta^{n-1}(a) - q^{n-1}
\delta^{n-1}(\sigma(a)) X$, we deduce from the induction hypothesis
that
$$h. \delta^n(a) = \lambda_h X \lambda_h^{n-1} \delta^{n-1}(h.a) -
q^{n-1} \lambda_h^{n-1} \delta^{n-1}(h.\sigma(a)) \lambda_h X.$$ Since
$h.\sigma(a)=hh_0.a=h_0h.a=\sigma(h.a)$, this leads to

\begin{eqnarray*}
h. \delta^n(a) &=& \lambda_h^n \left[ X \delta^{n-1}(h.a) - q^{n-1}
\delta^{n-1}(\sigma(h.a)) X \right]\\
& =& \lambda_h^n \left[ X
\delta^{n-1}(h.a) - \sigma \circ \delta^{n-1}(h.a) X \right]=
\lambda_h^n \delta^n(h.a).  
\end{eqnarray*}

This achieves the induction.\\$ $

Now, let  $a \in A$. Then, using the notations of \cite{c1}, we have
$$\theta (a) = \sum_{n=0}^{+ \infty} \frac{(1-q)^{-n}}{[n]!_q}
\delta^n \circ \sigma^{-n} (a) X^{-n}.$$

 Hence we get
$$h.\theta (a) = \sum_{n=0}^{+ \infty} \frac{(1-q)^{-n}}{[n]!_q}
h.\left[\delta^n \circ \sigma^{-n} (a)\right] h.X^{-n}.$$ Then the
previous study shows that
 $$h.\theta (a) = \sum_{n=0}^{+ \infty} \lambda_h^n
\frac{(1-q)^{-n}}{[n]!_q} \delta^n (h.\sigma^{-n} (a)) \lambda_h^{-n}
X^{-n}. $$ Now, since $\sigma$ coincide with the action of $h_0 \in \ch$
on $A$, we have
$h.\sigma^{-n}(a)=hh_0^{-n}.a=h_0^{-n}h.a=\sigma^{-n}(h.a)$, so that
$$h.\theta (a) = \sum_{n=0}^{+ \infty} \frac{(1-q)^{-n}}{[n]!_q}
\delta^n \circ \sigma^{-n} (h.a) X^{-n},$$ that is, $h.\theta (a) =
\theta(h.a)$ as desired.

\end{proof}

Note that for $b\in B$ with $b = \theta(a)$, we have $\alpha(b) =
\alpha(\theta(a)) = \theta(\sigma(a))= \theta(h_0.a) = h_0.\theta(a) =
h_0.b$; so that the action of $\alpha$ on $B$ coincides with the
action of $h_0$.

The above lemma shows that the action of $\ch$ on $\widehat{R}$
by automorphisms induces an action of $\ch$ on $B$ by
automorphisms. Further, since $T = B[X;\alpha]$ and since $X$ is an
$\ch$-eigenvector, this observation also proves that the action of $\ch$
on $\widehat{R}$ by automorphisms induces an action of $\ch$ on $T$ by
automorphisms.  Moreover, since every $\ch$-prime ideal of $A$ is
completely prime, we deduce that every
$\ch$-prime ideal of $B=\theta (A)$ is completely prime. Then, it
follows from  \cite[Proposition II.5.11]{bg} that every $\ch$-prime
ideal of $T = B[X;\alpha]$ is also completely prime. 

Let $b \in B$ be an $\ch$-eigenvector, say $h.b = \lambda_h b$ for
$\lambda_h \in k$, and suppose that $b = \theta(a)$. Then $\theta(h.a -
\lambda_h a) = h.\theta(a) -\lambda_h\theta(a) = h.b -\lambda_h b =
0$; so $h.a = \lambda_h a$ and $a$ is an $\ch$-eigenvector with the same
action of $\ch$ on $a$ as on $b$.

\begin{definition} \label{def-h-ufd}
{\rm Suppose that $A$ is a noetherian domain that is a $k$-algebra and
suppose that $\ch$ is a group acting on $A$ via
$k$-automorphisms. Then $A$ is an {\em $\ch$-UFD} if  each
nonzero $\ch$-prime $Q$ of $A$ contains a nonzero 
normal $\ch$-eigenvector $x$
such that the $\ch$-ideal $xA = Ax$ is completely prime.}
\sf{label:def-h-ufd}
\end{definition}

\begin{remark}{\rm  In particular, in an $\ch$-UFD, all $\ch$-primes of
height one as $\ch$-primes have height one as ordinary prime ideals, by the
principal ideal theorem. Thus, an ideal is an $\ch$-prime of height one as an
$\ch$-prime if and only if it is a prime $\ch$-ideal of height one as an
ordinary prime ideal. Also, in an $\ch$-UFD, the $\ch$-primes of height one
are principal, generated by a normal element, and completely prime.}
\end{remark}

\begin{proposition}
\label{proprec}
Let $R=A[X;\sigma,\delta]$ be a Cauchon extension. 
Suppose that $A$ is an $\ch$-UFD. Then $R$ is an $\ch$-UFD.
\end{proposition}

\begin{proof} Since $B$ is isomorphic to $A$ via $\theta$ and $\theta$ 
preserves the $\ch$-action, we know that every non-zero $\ch$-prime of $B$
contains a non-zero normal $\ch$-eigenvector $b$ such that $bB = Bb$ is a
completely prime ideal; that is, $B$ is an $\ch$-UFD. We start by showing that
such an element $b$ of $B$ can be used to produce, in a natural way, an
element of $R$ with similar properties.

Note, that every $\ch$-prime ideal of $A$ and $B$ is completely prime,
since this is one of the properties of $A$ being part of a Cauchon extension
and $B\cong A$ via a map compatible with the $\ch$-actions.

Let $b\in B$. Then $b\in B \sse T\sse \wr = RS^{-1}$; and so there
exists $n\geq 0$ with $bX^n \in R$.

Now, suppose that $ 0\neq b \in B$ is a normal $\ch$-eigenvector such
that $bB = Bb$ is a completely prime ideal. Choose $s\geq 0$ minimal
such that $x:= bX^s \in R$. We will show that $x$ is a normal
$\ch$-eigenvector in $R$ such that $xR = Rx$ is a completely prime
ideal.

First, note that $x$ is an $\ch$-eigenvector, since each of $b$ and $X$
is an $\ch$-eigenvector. Next, 
$$
Xb = \alpha(b)X = h_0.bX = \eta bX
$$ for some $0\neq \eta \in k$, since $b$ is an $\ch$-eigenvector. 

Hence, $b$ is normal in $T = B[X;\alpha]$. Also, $bT = Tb$ is an
$\ch$-invariant completely prime ideal of $T$.

It follows that $b\wr = \wr b$ is an $\ch$-invariant completely prime ideal of
$\wr$. However, $x\wr = bX^s\wr = b\wr$; and so $x\wr$ is an $\ch$-invariant
completely prime ideal of $\wr$. Thus, $I:= x\wr \cap R = b\wr \cap R$ is an
$\ch$-invariant completely prime ideal of $R$. We will show that $I = Rx$.

It is obvious that $Rx\sse I$. For the reverse inclusion, let $y\in I$. 
Then $y\in b\wr$ and so there exists $u\geq 0$ such that 
$yX^u \in bT = Tb$. 
Thus, 
there exists $c \in T$ such that $yX^u=cb$. 
Next, since $c
\in T \subseteq RS^{-1}$, there exists $v \geq 0$ such that
$cX^v \in R$. 
Set $r:=\eta^{-v} cX^v \in R$. 
Then, by using the fact that $Xb = \eta bX$, 
 we get
$yX^{u+v+s}=cbX^{v+s}=\eta^{-v} cX^{v}bX^{s}= r x$; and so there
exists $t \geq 0$ such that $y X^t = rx $ with $r\in R$. 
Choose such a $t$ minimal.

Assume that $t\geq 1$.
Express $r$, $y$ and $x$ as elements in the Ore
extension $R=A[X;\sigma,\delta]$, say,
$$r=\sum_{i=0}^d r_i X^i \mbox{, } \quad y=\sum_{i=0}^d y_i X^i \quad
\mbox{ and } \quad
x=\sum_{i=0}^d x_i X^i,
$$
where $d \geq 0$ and $r_i,y_i,x_i \in
A$ for all $ 0 \leq i\leq d $.
If $s=0$, then $x=b\in B\cap R \sse A$, and so $x_0 = b \neq 0$. If $s\geq 1$
then
$x_0 = 0$ would give
$$
bX^{s-1} = (bX^s)X^{-1} = xX^{-1} = \sum_{i=1}^d x_i X^{i-1} \in R,
$$
contradicting the minimality of $s$. Thus, $x_0 \neq 0$ whatever the value of
$s\geq 0$.

Recall that $Xb= \eta b X$, so that

$$
rx = \sum_{i=0}^d r_i X^i bX^s = \sum_{i=0}^d \eta^i r_i b X^{i+s}
= \sum_{i=0}^d \eta^i r_i x X^{i};
$$ 
that is,
$$
rx = \sum_{i,j=0}^d \eta^i r_ix_j X^{i+j}.
$$ 
Also, 
$rx=yX^t=\displaystyle{\sum_{i=0}^d y_i X^{i+t}}$; and so  we obtain 
the following equality
\begin{eqnarray}
\label{eqtoto}
\sum_{i,j=0}^d \eta^i r_ix_j X^{i+j} & = & \sum_{i=0}^d y_i X^{i+t}
\end{eqnarray}
in $R=A[X;\sigma,\delta]$.

Since $t \geq 1$, the term of degree 0 in the left hand side of
(\ref{eqtoto}) must be zero; that is, $r_0 x_0 =0$. 
Since $x_0
\neq 0$, this gives $r_0=0$. Hence $r=\sum_{i=1}^d r_i X^i= wX $
with $w=\sum_{i=1}^d r_i X^{i-1} \in R$. Consequently, the equality $y X^t =
r x$ can be rewritten as 
$$
y X^t = wX x=wXbX^s=\eta w bX^{s+1} = \eta w x X.
$$ 
It follows that 
$yX^{t-1} = \eta w x$, with $\eta w \in R$, contradicting 
the minimality of $t$.

Hence $t=0$ and  $y =rx$ with $r\in R$; so that $y \in Rx$, as required.

To sum up, we have established that $I=Rx$.\\

It remains to show that $xR = I$. First, note that, since $Xb=\eta b X$, we
have $x=\eta^{-s}X^s b \in R$ and it is clear that $\min\{i \in \mathbb{N}
\mid X^i b \in R \}=\min\{i \in \mathbb{N} \mid bX^i \in R \}=s$. Now by
writing elements of $R$ as polynomials with coefficients on the right, a very
similar calculation (which we omit) to that done above shows that $xR =I$.
Thus, $x = bX^s$ is a nonzero $\ch$-eigenvector of $R$ such that $I = xR = Rx
$ is a completely prime ideal. This finishes the first part of the proof.

Now, let $J$ be any nonzero $\ch$-prime ideal of $R$, and note that $J$ is
completely prime.

First, assume that $X\not\in J$. Then $JS^{-1} \cap T$ is a nonzero
$\ch$-invariant prime ideal of $T$ and it follows that $JS^{-1} \cap
B$ is a nonzero $\ch$-invariant prime ideal of $B$, by
Corollary~\ref{cont-ext-cor}.  Thus, there exists $0\neq b\in JS^{-1}
\cap B$ such that $b$ is a normal $\ch$-eigenvector and $bB = Bb$ is
a completely prime ideal of $B$. As in the earlier part of the proof, set
$x:= bX^s$, where $s$ is minimal such that $bX^s \in R$. Note that
$x\in JS^{-1} \cap R = J$, and that $x$ is a nonzero normal $\ch$-eigenvector
of $R$ such that $xR = Rx$ is a completely prime ideal of $R$.

Next, assume that $X\in J$. If $\delta = 0$ then $X$ is a nonzero normal
$\ch$-eigenvector such that $XR = RX$ is completely 
prime (since $A$ is a domain), as
required. Thus, we may assume that $\delta \neq 0$. 

Choose $c\in A$ such that $\delta(c) \neq 0$, and note that $0\neq \delta(c) =
Xc -\sigma(c) X \in J$; and so $J\cap A \neq 0$. It is clear that the map
$b\mapsto \theta^{-1}(b) +J$ defines a homomorphism from $B$ to $R/J$, and
this homomorphism extends to a homorphism $g$ from $T$ to $R/J$ such that
$g(X) =0$. This map, given by $g(\sum b_iX^i) = \theta^{-1}(b_0) +J$, 
commutes
with the action of $\ch$.
Set $J' = \ker(g)$; so that $J'$ is a completely prime $\ch$-ideal 
of $T$. With $c\in A$
as above, note that $g(\theta(\delta(c))) = \delta(c) + J = 0_{R/J}$. Thus,
$J'\cap B$ is a nonzero $\ch$-prime ideal of $B$. Thus, there is a nonzero
normal $\ch$-eigenvector $b\in J'\cap B$ such that $bB = Bb$ is 
a completely prime $\ch$-ideal 
of $B$. Set $x := bX^s$, where $s$ is minimal such that $bX^s \in R$.
Then, as in the earlier part of the proof, we know that $x$ is a nonzero
normal $\ch$-eigenvector of $R$ such that $xR = Rx$ is a prime ideal. In order
to finish this case, we will show that $x\in J$. Now, $b = \theta(a)$ for some
$0\neq a\in A$. We use the explicit formula for $\theta(a)$ to finish the
calculation:
\[b=\theta(a) = \sum_{n=0}^{+ \infty} \frac{(1-q)^{-n}}{[n]!_q} \delta^n
\circ \sigma^{-n} (a) X^{-n}.
\]
(The  sum on the right hand side exists since $\delta$ is
locally nilpotent). Since $\delta$ is locally nilpotent, there exists $d \in
\mathbb{N}$ such that $\delta^d(a) \neq 0$ and $\delta^{d+1}(a) = 0$. Then,
since $q \delta \circ \sigma = \sigma \circ \delta$, we have $$b=\theta(a) =
\sum_{n=0}^{d} \frac{(1-q)^{-n}}{[n]!_q} q^{n^2} \sigma^{-n} \circ \delta^n
(a) X^{-n},$$ and so the smallest integer $i$ such that $bX^i \in R$ is equal
to $d$. In other words, $s=d$ and $x=bX^d=\sum_{n=0}^{d}
\frac{(1-q)^{-n}}{[n]!_q} \delta^n \circ \sigma^{-n} (a) X^{d-n}$, that is: $$
x=\frac{(1-q)^{-d}}{[d]!_q} \delta^d \circ \sigma^{-d} (a) + \left(
\sum_{n=0}^{d-1} \frac{(1-q)^{-n}}{[n]!_q} \delta^n \circ \sigma^{-n} (a)
X^{d-1-n} \right)X.$$ Since $X \in J$, in order to prove that $x \in J$, it is
so sufficient to prove that $\delta^d \circ \sigma^{-d} (a)$ belongs to $J$.

Observe that, since $b \in J'$, we have $a=g(b) \in g(J') \subseteq J$. Hence,
if $d=0$, then $x=b=a$, and so $x \in J$ as desired. Assume now that $d \geq
1$. Then $\delta^d \circ \sigma^{-d} (a)=\delta \left( \delta ^{d-1} \circ
\sigma^{d}(a) \right)$. Set $e:=\delta ^{d-1} \circ \sigma^{d} (a) \in A$.
Then $\delta^d \circ \sigma^{-d} (a) = \delta(e) = Xe - \sigma(e)X \in J$,
since $X\in J$. This was what we needed to conclude that $x\in J$, as
required.

\end{proof}


\section{CGL extensions}

In this section, we develop a suitable context in which to apply the
results of the previous section to establish that certain iterated skew
polynomial extensions are $\ch$-UFDs.  The next problem is to use this
information, the Goodearl-Letzter stratification theory and the
noncommutative version of Nagata's lemma that we have established,
Proposition~\ref{pull-back-UFR}, to deduce that these extensions are,
infact, noetherian UFDs

The next definition contains all of the conditions that are necessary
for this programme to succeed. The definition is unwieldy, but is
justified by the fact that many of the quantum algebras that we wish
to study satisfy all of these conditions.

\begin{definition}\label{def-cgl}
{\rm  An iterated skew polynomial extension 
\[
A = k[x_1][x_2; \sigma_2, \delta_2] \dots [x_n; \sigma_n,\delta_n]
\]
is said to be a {\em CGL extension} (after Cauchon, Goodearl and Letzter)  
provided that the following list of
conditions is satisfied:

\begin{itemize}

\item With $A_j:= k[x_1][x_2; \sigma_2, \delta_2] \dots [x_j;
\sigma_j,\delta_j]$ for each $1\leq j\leq n$, each $\sigma_j$ is a
$k$-automorphism of $A_{j-1}$, each $\delta_{j}$ is a locally nilpotent
$k$-linear $\sigma_{j}$-derivation of
$A_{j-1}$, and there exist nonroots of unity $q_j \in k^*$ with
$\sigma_j\delta_j = q_j\delta_j\sigma_j$;

\item For each $i<j$ there exists a $\lambda_{ji}$ such that
$\sigma_j(x_i) = \lambda_{ji}x_i$;

\item There is a torus $\ch = (k^*)^r$ acting rationally on $A$ by $k$-algebra
automorphisms;

\item The $x_i$ for $1\leq i\leq n$ are $\ch$-eigenvectors;

\item There exist elements $h_1, \dots, h_n \in \ch$ such that $h_j(x_i) =
\sigma_j(x_i)$ for $j>i$ and such that the $h_j$-eigenvalue of $x_j$ is not a
root of unity. 

\end{itemize}

If, in addition, the subgroup of $k^*$ generated by the $\lambda_{ji}$ is
torsionfree then we will say that $A$ is a {\em torsionfree CGL extension}.}
\end{definition}

For a discussion of rational actions of tori, see \cite[Chapter II.2]{bg}.

Note that any CGL extension will be a noetherian domain with finite GK
dimension, cf. \cite[Lemma II.9.7]{bg}; 
and so will satisfy the descending chain condition on prime ideals,
as mentioned earlier. 
\\


Notice that, if $A$ is a CGL extension, then the action of ${\mathcal H}$ on
$k[x_1]$ is such that $h_1.x_1 =\lambda x_1$, where $\lambda\in k^\ast$ is not
a root of unity. From this, it follows easily that the only nonzero
${\mathcal H}$-prime of $k[x_1]$ is $\langle x_1\rangle$, which is
(completely) prime. Using [1, II.5.11], we deduce that, if $A$ is a CGL
extension 
then each of the extensions $A_j =
A_{j-1}[x_j;\sigma_j,\delta_j]$ is a Cauchon extension; so the results of the
previous section are available. Also, any CGL extension satisfies the
conditions of \cite[II.5.1]{bg} and so there are only finitely many
$\ch$-primes in $A$ and they are all completely prime, by \cite[Theorem
II.5.12]{bg}. Further, if $A$ is a torsionfree CGL extension, then all prime
ideals of $A$ are completely prime, by \cite[Theorem II.6.9]{bg}. In
particular, if such an $A$ is a noetherian UFR then it is a noetherian UFD.

\begin{proposition} \label{cgl-hufd}
Let $A$ be a CGL extension. Then $A$ is an $\ch$-UFD; that is, each nonzero
$\ch$-prime $Q$ of $A$ contains a nonzero 
normal $\ch$-eigenvector $a$ such
that the $\ch$-ideal $P:= aA = Aa$ is completely prime.
\end{proposition}

\begin{proof}
As already mentioned, the only nonzero
${\mathcal H}$-prime of $k[x_1]$ is $\langle x_1\rangle$ and it follows
immediately that $k[x_1]$ is an ${\mathcal H}$-UFD. Now,
each of the extensions $A_j = A_{j-1}[x_j;\sigma_j,\delta_j]$ is a
Cauchon extension; so apply Proposition~\ref{proprec} repeatedly.
\end{proof}

The main aim in this section is to show that any CGL extension is in fact a
noetherian UFR. It then follows that any torsionfree CGL extension is a
noetherian UFD. Since a CGL extension $A$ is an $\ch$-UFD, the prime ideals of
height one that are $\ch$-ideals are principal, generated by elements that are
normal and $\ch$-eigenvectors. Also, as noted above, there are only finitely
many $\ch$-primes, by \cite[Theorem II.5.12]{bg}, and they are all completely
prime. Thus, in order to show that such an extension is a noetherian UFD, we
have to deal with the primes of height one that are not $\ch$-primes. In the
language of Goodearl and Lezter, these primes are in the {\em stratum} of the
zero ideal; that is, if $P$ is a prime ideal of height one that is not an
$\ch$-prime, then the largest $\ch$-ideal contained in $P$ is the zero ideal.
The Goodearl-Letzter stratification theory enables us to deal with these
primes. The idea is simple. The stratification theory shows that, once we
invert all the regular $\ch$-eigenvectors, the prime ideals in the stratum of
the zero ideal become centrally generated. In fact, the height one primes in
the zero stratum become principal, generated by a central element in this
localisation; this shows this localisation is a noetherian UFR. However,
Proposition~\ref{pull-back-UFR} is valid only when we are inverting a
multiplicative set generated by finitely many normal elements. To deal with
this point, it turns out, and this is what we show first, that it is enough to
invert the multiplicative set generated by the finitely many generators of the
$\ch$-primes of height one in order to get a picture similar to that of the
stratification theory. \\

\begin{lemma} 
Let $I$ be an $\ch$-ideal in a CGL extension $A$. Then the prime
ideals minimal over $I$ are all $\ch$-prime ideals.
\end{lemma} 

\begin{proof}
Since
$A$ is noetherian, there are finitely many primes minimal over $I$. Let $Q$ be
a prime minimal over $I$. The ${\mathcal H}$-orbit of $Q$ consists of primes
minimal over $I$ and hence is finite. Now, \cite[II.2.9]{bg}
 shows that $Q$ is an
${\mathcal H}$-ideal.
\end{proof}

\begin{corollary}~~
Suppose that $A$ is a CGL extension and that $P_i = a_iA$ for $1\leq i
\leq m$ are the 
prime ideals of height one that are 
$\ch$-primes, where  
the $a_i$ are
normal $\ch$-eigenvectors. Then, 
each nonzero $\ch$-ideal of $A$
contains a product of the $a_i$ (repetitions allowed). 
\end{corollary}

\begin{proof} Let $I$ be a nonzero $\ch$-ideal of $A$.  Since $A$
is noetherian, there are only a finite number of prime ideals that are
minimal over $I$; denote these primes by $Q_1,\dots,Q_s$.  By the
previous lemma, these are all $\ch$-primes.  Since $A$ is noetherian,
the ideal $I$ contains a product of the $Q_i$. However, each $Q_i$
contains some $P_j$, by Proposition~\ref{cgl-hufd}; and so $I$
contains a product of the $P_i$, hence a product of the $a_i$. 
\end{proof}

Set $T$ to be the localisation of $A$ with respect to the
multiplicatively closed set generated by the normal $\ch$-eigenvectors
$a_i$. Then the rational action of $\ch$ on $A$ extends to an action of
$\ch$ on the localisation $T$ by $\k$-algebra automorphisms, since we
are localising with respect to $\ch$-eigenvectors, and this action of
$\ch$ on $T$ is also rational, by using \cite[II.2.7]{bg}. We have the
following proposition.

\begin{proposition}\label{t-is-hsimple}
~~The ring $T$ is $\ch$-simple; that
is, the only $\ch$-ideals of $T$ are $0$ and $T$. 
\end{proposition}

\begin{proof} Let $J$ be an $\ch$-ideal of $T$ and let $I=J \cap A$.
Clearly, $I$ is an $\ch$-ideal of $A$.
In addition, $J=IT$, by 
\cite[2.1.16]{mr}.  If $I=0$, then $J=0$. Otherwise, $J=T$, by
the previous corollary. 
\end{proof}

We are now in position to show that the CGL extension 
$A$ is a noetherian UFR.

\begin{theorem}\label{cgl-ufr}
~~Let $A = k[x_1][x_2; \sigma_2, \delta_2] \dots [x_n;
\sigma_n,\delta_n]$ be a CGL extension. Then $A$ is a noetherian UFR.  
\end{theorem}

\begin{proof} 
By  Proposition~\ref{pull-back-UFR}, it is enough to
prove that the localisation $T$ is a noetherian UFR. Now, as proved in
Proposition~\ref{t-is-hsimple}, $T$ is an $\ch$-simple ring. Thus, using
\cite[II.3.9]{bg}, it is a noetherian UFR, as required.
\end{proof}

\begin{theorem} \label{cgl-ufd}
Let $A$ be a torsionfree CGL-extension. Then 
$A$ is a noetherian UFD.
\end{theorem} 

\begin{proof}
Use Theorem~\ref{cgl-ufr} and the fact that all prime ideals are
completely prime in a torsionfree CGL-extension.
\end{proof}

This theorem applies to many quantum algebras. A selection of such algebras of
current interest is given in the following corollary. For exact definitions of
those of the algebras that are not explicitly defined in this paper, consult
\cite{g} or \cite[Section 6.2]{c1}

\begin{corollary} 
The following algebras are noetherian UFDs:
\begin{itemize}

\item The algebra of quantum matrices $\oqmmn$, with $q$ not a root of unity,
(see also the next section for more information about $\oqmmn$), and, more
generally, the multiparameter version $\co_{\lambda, {\bf p}}(M_{m,n}(k))$,
with 
$\lambda$ not a root of unity and the group 
$\langle {\lambda, p_{ij}}\rangle$ torsionfree.

\item The quantized enveloping algebra $U_{q}({\bf \frak{n}}^+)$, with $q$
not a root of unity, 
of the nilpotent subalgebra ${\bf \frak{n}}^+$ of a complex
semisimple Lie algebra $\frak{g}$.

\item The quantized enveloping algebra $U_{q}({\bf \frak{b}}^+)$, with $q$
not a root of unity, 
of the Borel subalgebra ${\bf \frak{b}}^+$ of a complex
semisimple Lie algebra $\frak{g}$.

\item The quantum affine space $\co_{\bf q}(k^n)$, with 
$\langle {q_{ij}}\rangle$ torsionfree.



\item The quantized Weyl algebra $A_n^{Q,\Gamma}(k)$ with each $q_i$ not a
root of unity and $\langle {q_i, \gamma_{ij}}\rangle$ torsionfree.

\item The quantum grassmannian $\gqmn$, with $q$ not a root of unity.

\end{itemize} 

\end{corollary} 

\begin{proof} 
The algebras $\oqmmn, \co_{\lambda, {\bf p}}(M_{m,n}(k)), \co_{\bf q}(k^n), 
A_n^{Q,\Gamma}(k)$ are described in \cite{g} as iterated skew polynomial
extensions with appropriate torus actions, and can easily be checked to be
torsionfree CGL-extensions. (The only awkward point is to check that the first
condition holds, and, in particular, to check that the $\delta_i$ involved all
act locally nilpotently. The lemma below, which is easy to prove,  
helps deal with this point.)

The algebra $U_{q}({\bf \frak{n}}^+)$ is described in \cite[Section 6.2]{c1}
and is easily seen to be a CGL-extension.

The algebra $U_{q}({\bf \frak{b}}^+)$ is described in \cite{g} as a
localisation of an algebra that is an iterated skew polynomial extension with
a torus action. This algebra is easily checked to be a CGL-extension.

The algebra $\gqmn$ is shown to be a noetherian UFD in the last section of
this paper.

\end{proof}

\begin{lemma} 
Let $R$ be a $k$-algebra,
$\tau$ a $k$-algebra automorphism, $\delta$ a left $\tau$-derivation,
which we assume to be $k$-linear and set $S=R[x;\tau,\delta]$. In 
addition, let $X \subseteq R$ be a generating set of the $k$-algebra $R$. 
Then, the following holds.\\
(i) Assume that there exists $q \in k$ such that, for all $x \in X$,
$\delta\tau(x)=q\tau\delta(x)$, then $\delta\tau=q\tau\delta$.\\
(ii) Assume that there exists $q \in k$ such that 
$\delta\tau=q\tau\delta$.
If, for all $x \in X$, there exists $d\in{\mathbb N}^\ast$ such that 
$\delta^{d}(x)=0$,
then $\delta$ is locally nilpotent.
\end{lemma}

\section{Height one $\ch$-primes in $\oqmmn$}


In this section, we identify generators for each of the 
height one primes which are ${\mathcal H}$-ideals of 
the algebra of quantum matrices, in the generic
case.

Throughout, $\k$ is a field and $q$ is a nonzero element of $\k$ that
is not a root of unity.  Let $m,n$ be positive integers.  Recall that
the algebra of $m\times n$ quantum matrices, $\oqmmn$, is the
$\k$-algebra generated by $mn$ indeterminates $x_{ij}$, with $1 \leq i
\leq m$ and $1\leq j\leq n$, subject to the relations
\[
\begin{array}{ll}
x_{ij}x_{il}=qx_{il}x_{ij},
& (j < l); \\
x_{ij}x_{kj}=qx_{kj}x_{ij},
& (i<k); \\
x_{ij}x_{kl}=x_{kl}x_{ij},
& (i <k,  j>l); \\
x_{ij}x_{kl}-x_{kl}x_{ij}=(q-q^{-1})x_{il}x_{kj},
& (i<k,  j<l). 
\end{array}
\]
In the case that $m=n$, we write $\oqmn$ for $\oqmmn$.  

In view of the restriction that $q$ is not a root of unity, we refer to
$\oqmmn$ as the algebra of {\em generic quantum matrices}.

Let $\ch$ be the $(m+n)$-torus $(k^*)^m \times (k^*)^n$. The torus $\ch$
acts on $\oqmmn$ by $k$-algebra automorphisms in the following way:
\[
(\alpha_1, \dots , \alpha_m, \beta_1, \dots, \beta_n)\cdot x_{ij}:=
\alpha_i \beta_j x_{ij}.
\]

The algebra $\oqmmn$ can be presented as an iterated skew polynomial
extension with the variables added in lexicographical order. With this
presentation, and with the group $\ch$ above acting, $\oqmmn$ is a
torsionfree CGL extension; and so is a noetherian UFD by the results of the
previous section. There are only finitely many 
height one prime ideals which  are $\ch$-primes,  
and the purpose of this section is to identify these
$\ch$-primes.

In the literature, many results are only 
stated for $\oqmn$ but are easily translated
to $\oqmmn$, by using arguments based on the following easy
observations. First, if $I$ is a set of row indices and $J$ is a set of
column indices then the subalgebra of $\oqmn$ or $\oqmmn$ generated by
the $x_{ij}$ with $i\in I$ and $j\in J$ is isomorphic to another quantum
matrix algebra in a natural way. Secondly,  let $A=\oqmn$, and let
$B=\oqmmn$, with $m\leq n$, 
be the quantum matrix algebra generated by generators in the 
first $m$ rows  of $A$, then there is an algebra epimorphism 
$\pi \, : \, A \longrightarrow B$ defined by the projection
given by $x_{ij} \mapsto x_{ij}$ if $i \le m$ and $x_{ij} \mapsto 0$
otherwise. 
By using the first observation, we may think of $\oqmmn$ and $\oqmnm$
being embedded in a common $\oqmn$. Then, there is an
isomorphism between $\oqmmn$ and $\oqmnm$ given by transposition of the
generators in $\oqmn$; that is, $x_{ij}\mapsto x_{ji}$, see \cite[Proposition 
3.7.1]{pw}.  For this reason, we will assume that 
$m\leq n$.  In view of the restriction that $q$ is not a root of unity,
we will refer to $\oqmmn$ as a {\em generic quantum matrix algebra}.

The algebra $\oqmn$ has a special element, $\detq$, the {\em quantum
determinant}, defined by 
\[
\detq := \sum_{\sigma}\, (-q)^{l(\sigma)}x_{1\sigma(1)}\cdots
x_{n\sigma(n)}, 
\]
where the sum is taken over the permutations of $\{1, \dots, n\}$ and
$l(\sigma)$ is the usual length function on such permutations. The
quantum determinant is a central element of $\oqmn$, see, for example,
\cite[Theorem 4.6.1]{pw}. 
If $I$ is a $t$-element subset of $\{1, \dots, m\}$ and
$J$ is a $t$-element subset of $\{1, \dots, n\}$ then the 
quantum determinant of the subalgebra of
$\oqmmn$ generated by $\{x_{ij}\}$, with $i\in I$ and $j\in J$, is
denoted by 
$[I\mid J]$. The elements $[I\mid J]$ are the {\em
quantum minors} of $\oqmmn$.  They are not in general central; however,
they do possess good commutation properties: in particular, in what
follows, we will 
identify several quantum minors that are normal elements. Two elements
$a, b$ are said to {\em $q$-commute} if there is an integer $s$ such
that $ab = q^sba$. An element that $q$-commutes with each of the
generators of a quantum matrix algebra is easily seen to be normal, and
this is a standard way to demonstrate normality. In many sources, such
commutation relations are established for $\oqmn$. Usually, it is easy
to transfer such results to $\oqmmn$, by including this quantum matrix
algebra as a subalgebra of a suitable $\oqmn$ by including extra rows or
columns of generators: obviously, if an element $q$-commutes with each
of the generators in this larger algebra then it $q$-commutes with the
generators of the original algebra. In addition, we will use the
transposition isomorphism to derive further $q$-commutation results,
with little comment.

Cauchon's theory of deleting derivations,\cite{c1, c2}, has been
applied to quantum matrices with great success.
In fact, in \cite{c2}, Cauchon works with $\oqmn$;
however, the methods extend to $\oqmmn$ and the details are worked out
in \cite{laun}. Let $\w$ denote an $m\times n$ array of square boxes in
which each box is coloured either black or white. A {\em Cauchon
diagram} is such an array with the following property: if a square is
coloured black then either every square to the left of this square is
also coloured black, or every square above this square is also coloured
black. Cauchon \cite{c2} and Launois \cite{laun}  prove that the
$\ch$-prime ideals of $\oqmmn$ are in bijection with the $m\times n$
Cauchon diagrams. In addition, if $P$ is an $\ch$-prime, then the height
of $P$ (as a prime ideal) is equal to the number of black boxes in the
corresponding diagram, by \cite{c2},~Th\'eor\`eme~6.3.3 (which is easily
adapted to the rectangular case), and \cite{laun},~Proposition~1.3.2.2. 
(Recall that, by [1, II.2.9], any ${\mathcal H}$-prime is prime.)\\
 
For $1\leq i \leq m$, let $c_i$ denote the $i\times i$ quantum minor $[m-i+1,
\dots , m\mid 1, \dots, i]$ of $\oqmmn$ and let $b_i$ denote the $i\times i$
quantum minor $[1, \dots , i\mid n -i+ 1, \dots, n]$ of $\oqmmn$, while for
$m< i \leq n$, let $b_i$ denote the quantum minor $[1, \dots, m \mid n-i+1,
\dots, n+m -i]$. Note that $c_m = b_n$; in particular, note that in $\oqmn$ we
have $b_n = c_n = \detq$. For orientation, note that the $b_i$ are the minors
coming from the top right of the matrix of generators of $\oqmmn$, while the
$c_i$ come from the bottom left.\\

The quantum minors are $\ch$-eigenvectors; and so, for example, the ideals
generated by each of the elements $b_i$ and $c_i$, defined above, are
$\ch$-ideals. We will show below that they are $\ch$-prime ideals. 

\begin{lemma}\label{bcnormal}
The elements $b_i$, with $1\leq i\leq n$, and $c_i$,  with $1\leq i\leq
m$,  are normal elements of $\oqmmn$. 
\end{lemma}

\begin{proof}
Let $1\leq i\leq m$; it follows easily from \cite[Corollary
A.2]{glen} that $c_i$ $q$-commute with each generator of $\oqmmn$ and,
using the transpose automorphism, that the same is true for $b_i$. For
$b_i$, with $m < i\leq n$, a slightly more complicated argument is
required. Fix an $i$ with $m < i\leq n$. Consider a generator
$x_{kl}$. If $l\leq n+m-i$ then $x_{kl}$ and $b_i$ belong to the
quantum matrix algebra $\oq(M_{m,n+m-i})$ obtained from the generators
in the first $n+m-i$ columns of $\oqmmn$. In fact, $b_i$ is $b_m$ in
this subalgebra, and so $x_{kl}$ and $b_i$ $q$-commute. If $l>n+m-i$
then $x_{kl}$ and $b_i$ belong to the quantum matrix algebra
$\oq(M_{m,i})$ obtained from the generators in the last $i$ columns of
$\oqmmn$. In this case, $b_i$ is $c_m$ in this subalgebra and so again
we see that $x_{kl}$ and $b_i$ $q$-commute. Thus $b_i$ $q$-commutes
with each of the generators of $\oqmmn$ and so is a normal element in
this algebra.
\end{proof}

\begin{proposition}\label{m+n-1}
~~There are precisely $m+n-1$ height one primes
that are $\ch$-primes in the generic quantum matrix algebra $\oqmmn$. 
They are the ideals generated by $b_1,\dots, b_n$
and $c_1, \dots, c_{m-1}$ (recall that $c_m = b_n$). 
\end{proposition}

\begin{proof}~~It is easily seen that the elements $b_1,\dots, b_n,
c_1, \dots, c_{m-1}$ generate pairwise distinct ideals.

The
height one primes that are ${\mathcal H}$-primes
are in bijection with the Cauchon diagrams with precisely one
black box. Such Cauchon diagrams arise by filling in one box either in
the first row of the array, or the first column. There are $m+n-1$ ways
of doing this; and so there are $m+n-1$ 
height one primes that are $\ch$-primes.

That the ideals specified are $\ch$-ideals is due to the fact that the
$b_i$ and $c_i$ are $\ch$-eigenvectors.  That the ideals are prime comes
about in the following way.  If we restrict to the quantum submatrix
algebra $A$, say, specified by the rows and columns of a $b_i$ or $c_i$,
then when we factor out $b_i$ or $c_i$ from $A$ we are factoring out the
quantum determinant of $A$, and so the factor $A/b_iA$ or $A/c_iA$ is a
domain, see, for example, \cite[Theorem 2.5]{glen}.  Since the $b_i$ or
$c_i$ $q$-commute with the remaining $x_{ij}$ we can add the remaining
$x_{ij}$ in such a way that at any stage if we have reached a subalgebra
$B$ then $B/b_iB$, say, is an iterated skew polynomial algebra over
$A/b_iA$ and so is a domain.  For example, if we are in the case that
$m<i\leq n$, then we can add the $x_{ij}$ to the left of the rows and
columns used by $b_i$ by moving from right to left along each row,
starting with the bottom row and moving upwards row by row.
We then can add the $x_{ij}$ 
to the right of the rows and columns used by $b_i$ in lexicographic
order.  Thus each $\oqmmn/b_i\oqmmn$ and $\oqmmn/c_i\oqmmn$ is a domain
and so each ideal of $\oqmmn$ generated by a $b_i$ or $c_i$ is a
completely prime ideal.  Since these ideals are $\ch$-ideals, they are
also $\ch$-primes.  Since we have precisely $m+n-1$ elements $b_i$ or
$c_i$ this gives all of the 
height one primes that are $\ch$-primes.
\end{proof}

\section{Generic quantum grassmannians are UFD}

Recall that the {\em quantum grassmannian} subalgebra, $\gqmn$, of $\oqmmn$ is
the subalgebra generated by the $m\times m$ maximal quantum minors of $\oqmmn$
(recall that we are assuming that $m\leq n$). The algebra $\gqmn$ is a
noetherian domain, see, for example, \cite[Theorem 1.1]{klr}. Our usual
restriction that $q$ is not a root of unity applies in this section; so we
refer to $\gqmn$ at the {\em generic quantum grassmannian}.

In view of the fact that each of the quantum minors that generates 
$\gqmn$ is of the form $[1, \dots, m \mid J]$ we will denote such a
minor by $[J]$. The two extreme  quantum minors, $[1, \dots, m]$ and
$[n-m+1, \dots, n]$ are normal in $\gqmn$, see, for example,
\cite[Corollary 1.1, Lemma 1.1]{klr}. 

We will use the fact that generic quantum matrices are noetherian UFD, and the
dehomogenisation isomorphism 
\[ \oq(M_{m,n-m})[y,y^{-1};\phi] \goesto \gqmn
[[n-m+1, \dots, n]^{-1}] 
\] 
of \cite[Corollary 4.1]{klr} to show that $\gqmn$ is a noetherian UFD. Note
that the automorphism $\phi$ used in the dehomogenisation isomorphism acts on
generators via $\phi(x_{ij}) = q^{-1}x_{ij}$, see \cite[Corollary 4.1]{klr}.

To show that $\gqmn$ is a noetherian UFD, we proceed as follows. 
First, we show
that the problem reduces to proving that the localisation 
$\gqmn [[n-m+1, \dots, n]^{-1}]$ 
is a noetherian UFD. Once this is done, by the
dehomogenisation theorem, the problem transfers to showing that 
$\oq(M_{m,n-m})[y,y^{-1};\phi]$ is a noetherian UFD, and this is the
second step.\\  

The first step is easy and essentially amounts to proving Lemma \ref{u-is-cp}
below.
 
\begin{lemma}\label{u-is-cp}The ideal of $\gqmn$ generated by $[n-m+1,
\dots, n]$ is a completely prime ideal.  
\end{lemma} 

\begin{proof}~~Let $R = \gqmn$. 
The isomorphism discussed immediately before \cite[Lemma
1.1]{klr} shows that the result follows provided that we show that $a:=[1,
\dots, m]$ generates a completely prime ideal in $\gqmn$. 

Note that $u:= [n-m+1,\dots, n]$ is left regular modulo $aR$, see
the proof of Theorem 6.1 of \cite{klr}.  Hence, it is enough to prove
that $a$ generates a completely prime ideal in the localisation
$R[u^{-1}]$.  We use the dehomogenisation isomorphism introduced above. 
Set $v:= [1,\dots,t|1,\dots,t]$ with $t=m$ if $m \le n-m$ and $t=n-m$
otherwise. 
By \cite[Lemma 3.5.1]{lr}, it is enough to show that $v$ 
generates a completely prime ideal of
$\oq(M_{m,n-m})[y,y^{-1};\phi]$.  However, $v$ generates a completely
prime ideal of $\oq(M_{m,n-m})$ that is left invariant by $\phi$, since
$\phi(v)$ is a scalar multiple of $v$.  Thus $v$ generates a completely
prime ideal of $\oq(M_{m,n-m})[y,y^{-1};\phi]$, as required. 
\end{proof} 

To achieve the second step, we observe first that $\oq(M_{m,n-m})[y;\phi]$ is
a torsionfree CGL-extension (in $s:= m(n-m)\!+\!1$ steps) as follows. The
torus $\ch = (k^*)^n = (k^*)^m \times (k^*)^{n-m}$ acts on $\oq(M_{m,n-m})$ as
defined at the beginning of the previous section, and we have already observed
that this makes $\oq(M_{m,n-m})$ a CGL-extension (in $s-1 = m(n-m)$ steps). In
order to deal with the last step (extension by $y$) we proceed as follows. We
extend this action of $\ch$ to $\oq(M_{m,n-m})[y;\phi]$ by setting
$(a_1,\dots,a_n).y=a_1 \dots a_n y$. The element $h_{s}$ needed for the final
extension is given by $h_{s}:=(q^{-1},\dots,q^{-1},1,\dots,1) \in \ch$ (with
$m$ occurences of $q^{-1}$), since we require that $h_{s}(x_{ij}) =
\phi(x_{ij}) = q^{-1}x_{ij}$. Moreover $h_{s}.y=q^{-m} y$, and $q^{-m}$ is not
a root of unity, since $q$ is not. With this information provided, it is easy
to check the remaining conditions and conclude that $\oq(M_{m,n-m})[y;\phi]$
is a torsionfree CGL-extension.

\begin{theorem}\label{laurent-is-factorial}
Suppose that $q \in {\k^\ast}$ is not a
root of unity. Then $\oq(M_{m,n-m})[y,y^{-1};\phi]$ is
a noetherian UFD. 
\end{theorem}

\begin{proof}
That $\oq(M_{m,n-m})[y; \phi]$ is a noetherian UFD 
follows from
Theorem~\ref{cgl-ufd}, since $\oq(M_{m,n-m})[y; \phi]$ is a torsionfree 
CGL-extension.
It follows that $\oq(M_{m,n-m})[y, y^{-1}; \phi]$ is a
noetherian UFD.
\end{proof}

\begin{theorem}\label{g-is-factorial}
The generic quantum grassmannian, $\gqmn$, is a noetherian UFD. 
\end{theorem}

\begin{proof}~~The previous result shows that $\gqmn[[n-m+1, \dots,
n]^{-1}]$ is a noetherian UFD, by using the dehomogenisation isomorphism.  By
Lemma~\ref{u-is-cp}, we know that $[n-m+1, \dots, n]$ generates a
completely prime ideal.  Thus the result follows from
Lemma~\ref{nagata}.  
\end{proof}

\noindent{\bf Acknowledgment}~We thank Gerard Cauchon for many useful
comments.



\vskip 1cm

\newpage

\noindent S Launois:\\
Laboratoire de Math\'ematiques - UMR6056, Universit\'e de Reims\\
Moulin de la Housse\\  BP 1039\\ 51687 REIMS Cedex 2\\ France\\
E-mail : stephane.launois@univ-reims.fr
\\

\noindent T H Lenagan: \\
School of Mathematics, University of Edinburgh,\\
James Clerk Maxwell Building, King's Buildings, Mayfield Road,\\
Edinburgh EH9 3JZ, Scotland\\
E-mail: tom@maths.ed.ac.uk \\
\\
L Rigal: \\
Universit\'e Jean Monnet
(Saint-\'Etienne), \\
Facult\'e des Sciences et
Techniques, \\
D\'e\-par\-te\-ment de Math\'ematiques,\\
23 rue du Docteur Paul Michelon,\\
42023 Saint-\'Etienne C\'edex 2,\\France\\
E-mail: Laurent.Rigal@univ-st-etienne.fr

\end{document}